\documentclass[reqno]{amsart}

\usepackage[utf8]{inputenc}
\usepackage[T1]{fontenc}
\usepackage[american]{babel,isodate}
\cleanlookdateon
\usepackage{lmodern}
\usepackage{hyperref}
\usepackage{doi}
\usepackage{xcolor}
\usepackage{stmaryrd}
\usepackage{enumerate}
\usepackage{listings}
\usepackage{amsmath,amssymb,xcolor}
\definecolor{lightgrey}{rgb}{0.9,0.9,0.9}
\usepackage{physics}
\usepackage{ stmaryrd }
\usepackage{tikz}
\usepackage{tikz-cd}
\usepackage{wrapfig}
\usetikzlibrary{arrows}
\usepackage{cancel}
\usetikzlibrary{babel}
\usepackage{lscape}
\usepackage{url}

\usepackage{hyperref}
\hypersetup{linktocpage,colorlinks=true,citecolor=purple,linkcolor=blue,urlcolor=blue,filecolor=black}

\usepackage{float}
\usepackage{graphicx, import}
\usepackage{mathtools}
\usepackage{dsfont}
\usepackage{pinlabel}
\usepackage{yhmath}

\usepackage{xargs}                      
\usepackage{xcolor} 
\definecolor{OliveGreen}{rgb}{0,0.6,0}
 
\usepackage[colorinlistoftodos, prependcaption,textsize=tiny]{todonotes}
\newcommandx{\unsure}[2][1=]{\todo[linecolor=red,backgroundcolor=red!25,bordercolor=red,#1]{#2}}
\newcommandx{\change}[2][1=]{\todo[linecolor=blue,backgroundcolor=blue!25,bordercolor=blue,#1]{#2}}
\newcommandx{\info}[2][1=]{\todo[linecolor=OliveGreen,backgroundcolor=OliveGreen!25,bordercolor=OliveGreen,#1]{#2}}

\newcommand{\centre}[1]{\begin{array}{c} #1 \end{array}}

\usepackage{stackengine}

\usepackage{amsmath,amsthm,amssymb, amsfonts, amscd}
\usepackage{bbm}
\usepackage{rotating,subcaption}
\usepackage{url}
\usepackage{graphicx,bm}
\usepackage{mathtools}
\usepackage[mathscr]{euscript}
\allowdisplaybreaks
\usepackage{lipsum}
\usepackage{bm}
\usepackage[capitalize]{cleveref}

\newcommand{\Det}{\operatorname{Det}}

\usepackage{xpatch}
\makeatletter
\AtBeginDocument{\xpatchcmd{\@thm}{\thm@headpunct{.}}{\thm@headpunct{}}{}{}}
\makeatother

\newtheorem{theorem}{Theorem}[section]
\newtheorem{proposition}[theorem]{Proposition}
\newtheorem{propdef}[theorem]{Proposition/Definition}
\newtheorem{lemma}[theorem]{Lemma}

\theoremstyle{definition}
\newtheorem{definition}[theorem]{Definition}
\newtheorem{example}[theorem]{Example}

\theoremstyle{remark}

\newtheorem{remark}[theorem]{Remark}
\newtheorem{notation}[theorem]{Notation}

\numberwithin{equation}{section}

\newcommand{\uno}{\mathds{1}}
\newcommand{\A}{\mathcal{A}}

\newcommand{\bk}{\Bbbk}

\usepackage{accents}
\newcommand{\uhat}{\underaccent{\check}}

\makeatletter
\newcommand{\cupr@tip}{\text{\raisebox{-0.1ex}{$\m@th\hat{}$}}}
\newcommand{\cupr}{\mathbin{\cup\cupr@}}

\newcommand{\cupr@}{%
  \mathchoice
  {\mkern-1.35mu\cupr@tip}
  {\mkern-1.35mu\cupr@tip}
  {\mkern-1.55mu\cupr@tip}
  {\mkern-1.875mu\cupr@tip}
}
\makeatother

\makeatletter
\newcommand{\capr@tip}{\text{\raisebox{0.47ex}{$\m@th\uhat{}$}}}
\newcommand{\capr}{\mathbin{\capr@\cap}}

\newcommand{\capr@}{%
  \mathchoice
  {\mkern11.6mu\capr@tip\mkern-11.6mu}
  {\mkern11.4mu\capr@tip\mkern-11.4mu}
  {\mkern11.1mu\capr@tip\mkern-11.1mu}
  {\mkern10.2mu\capr@tip\mkern-10.2mu}
}
\makeatother

\makeatletter
\newcommand{\capl@tip}{\text{\raisebox{0.47ex}{$\m@th\uhat{}$}}}
\newcommand{\capl}{\mathbin{\capl@\cap}}

\newcommand{\capl@}{%
  \mathchoice
  {\mkern2.1mu\capl@tip\mkern-2.1mu}
  {\mkern2.1mu\capl@tip\mkern-2.1mu}
  {\mkern2.3mu\capl@tip\mkern-2.3mu}
  {\mkern2.1mu\capl@tip\mkern-2.1mu}
}
\makeatother

\makeatletter
\newcommand{\cupl@tip}{\text{\raisebox{-0.1ex}{$\m@th\hat{}$}}}
\newcommand{\cupl}{\mathbin{\cupl@\cup}}
\newcommand{\cupl@}{%
  \mathchoice
  {\mkern1.35mu\cupl@tip\mkern-1.35mu}
  {\mkern1.35mu\cupl@tip\mkern-1.35mu}
  {\mkern1.55mu\cupl@tip\mkern-1.55mu}
  {\mkern1.875mu\cupl@tip\mkern-1.875mu}
}
\makeatother

\DeclareFontFamily{U}{mathx}{}
\DeclareFontShape{U}{mathx}{m}{n}{ <-> mathx10 }{}
\DeclareSymbolFont{mathx}{U}{mathx}{m}{n}
\DeclareFontSubstitution{U}{mathx}{m}{n}
\DeclareMathAccent{\widecheck}{0}{mathx}{"71}

\begin{document}


\title{Tensors, Gaussians and the Alexander Polynomial}

\date{\today}

\author{Boudewijn Bosch}
\address{Bernouilli Institute, University of Groningen, Nijenborgh 9, 9747 AG, Groningen, The Netherlands}
\email{\href{mailto:b.j.bosch@rug.nl}{b.j.bosch@rug.nl}}




\begin{abstract}
 In this paper, we formally show that the Alexander polynomial of a knot can be computed through contractions of a Gaussian function. We build on the approach of Bar-Natan and Van der Veen to universal knot invariants using (perturbed) Gaussian functions. We endow the Weyl--Heisenberg algebra with an XC-structure. Given a knot, the universal invariant arising from the XC-algebra induces a Gaussian function whose partition function recovers the Alexander polynomial. This formalism allows for the addition of suitable perturbations to the Gaussian function, from which perturbative expansions of knot invariants can be derived.

\end{abstract}

\subjclass{57K14, 16T05, 17B37}


\maketitle



\section{Introduction}
\noindent
The Alexander polynomial is a classical invariant of oriented knots in $S^3$ \cite{Alexander1928TopologicalIO}. To each knot $\mathcal K$ it assigns a Laurent polynomial $\Delta_{\mathcal K}(T) \in \mathbb Z[T^{\pm 1}]$. It can be defined in many equivalent ways, such as via covering spaces, Seifert matrices, Fox calculus, etc. In this chapter we will rely on the description coming from a presentation matrix of the Alexander module of a knot; see, e.g., \cite[Ch. 8]{rolfsen2003knots}.

Let
\[
A \in \operatorname{Mat}_{N}(\mathbb{Z}[T^{\pm 1}])
\]
be a specifically chosen presentation matrix of the Alexander module for a knot $\mathcal K$ (see Eq. \eqref{eq:alexandermatrix}). The Alexander polynomial is retrieved as $\Delta_{\mathcal K}(T) = \Det{A}$, up to multiplication by $\pm T^n$. In \cite{barnatan2024perturbedalexanderinvariant} a “perturbed-Alexander invariant” is defined as a quadratic expression in the entries of $A^{-1}$. The construction in \cite{barnatan2024perturbedalexanderinvariant} implicitly used a ``perturbed'' $R$-matrix, which, according to the authors, gives rise to formulas for the first-order perturbed Alexander invariant.

Notably, these formulas can be interpreted as arising from a Gaussian integral.  The matrix $A(T)$ plays the role of a precision matrix of a Gaussian function, and perturbations can be derived using Wick's theorem. However, this point of view is also only implicitly used in \cite{barnatan2024perturbedalexanderinvariant}.

The purpose of this chapter is to make this Gaussian picture precise and to use it as a basis for constructing the Alexander polynomial. We build on the approach of Bar-Natan and Van der Veen to universal knot invariants through perturbed-Gaussian functions \cite{bar2021perturbed}. By considering the universal knot invariant defined through XC-algebras, as introduced in \cite{becerra2025refinedfunctorialuniversaltangle}, together with the Weyl--Heisenberg algebra, we construct a Gaussian function whose partition function is shown to recover the Alexander polynomial $\Delta_{\mathcal K}(T)$ (see Theorem \ref{thm:alex}). The use of the Heisenberg algebra was originally introduced in \cite[Section 3.1]{bar2021perturbed}, where its connection to the Alexander polynomial was suggested through examples.

We introduce two XC-algebras, denoted by $\mathcal{A}$ and $\mathcal{B}$, whose universal $R$-matrix is related through the notion of ``twisting'' as introduced in  \cite{bosch2025largecolorexpansionderiveduniversal}. The XC-algebra $\mathcal{B}$ is recognized from the perspective of the Burau representation of braid group generators, while the XC-algebra $\mathcal{A}$ has origins within the representation theory of quantum groups and quantum cluster algebras.

A benefit of a Gaussian approach to the Alexander polynomial is that it allows for the derivation of perturbed Alexander invariants by adding suitable perturbations to the Gaussian function. This is expanded on in Section \ref{sec:discussion}.

\subsection{Organization of the paper}
Section~\ref{sec:tenscont} introduces tensor products and multiplication maps. In Section~\ref{sec:heisenberg} we develop the Heisenberg algebra framework and set up the normal-ordering and contraction formalism. In Section~\ref{sec:alexander} we construct an XC-algebra from which we derive the Alexander polynomial. Finally, in Section~\ref{sec:discussion} we discuss possible future directions of research.

\subsection*{Acknowledgments} The author would like to thank Roland van der Veen for many helpful discussions and suggestions on the content of this paper.

\section{Tensor Products}
\label{sec:tenscont}

\noindent
In this section, we introduce notation we use throughout this chapter. We assume that the reader is familiar with tensor products, but may not have seen them in this particular notation.

Let $\bk$ be a commutative ring and let $A$ be a unital associative $\bk$-algebra. For any finite set $I$ we write
\[
A^{\otimes I} := \bigotimes_{i\in I} A,
\]
Note that any bijection $\sigma\colon I\xrightarrow{\sim} J$ induces an algebra isomorphism
\[
P_\sigma\colon A^{\otimes I}\xrightarrow{\ \cong\ }A^{\otimes J},
\qquad
P_\sigma\left(\bigotimes_{i\in I} a_i\right)
=
\bigotimes_{j\in J} a_{\sigma^{-1}(j)}.
\]

\begin{definition}\label{def:leg}
For $i\in I$, the \emph{leg embedding}
\[
\iota_i^I\colon A\longrightarrow A^{\otimes I}
\]
is defined as the unital algebra map which sends $a\mapsto \bigotimes_{j\in I} x_j$ with $x_i=a$ and $x_{j}=1$ for $j \neq i$. \end{definition}
\begin{notation}
We denote
\[
a_i := \iota_i^I(a) \in A^{\otimes I}.
\]
Moreover, for a finite subset $S\subset I$ and $\{a_s\in A\}_{s\in S}$ we write $\prod_{s\in S} a_s$ as the product in $A^{\otimes I}$.
\end{notation}

The following is obvious.
\begin{proposition}\label{prop:leg-relations}
For all $a,b\in A$ and all $i\neq j$ in $I$ one has
\[
a_ib_j=b_ja_i,\qquad a_ib_i=(ab)_i .
\]
Consequently, $\prod_{s\in S}a_s$ is independent of the order of multiplication.
\end{proposition}

\begin{example}\label{ex:two-legs}
Let $I = \{i,j\}$. For $a,b,c\in A$ one has
\[
a_i b_j c_i=(ac)_ib_j,\qquad
(a+b)_i c_j=a_i c_j + b_i c_j,
\]
and the element $a_i b_j$ amounts to the regular tensor product $a\otimes b$.
\end{example}

We next define the subspace on which a choice of leg is the unit.

\begin{definition}\label{def:Vk}
Let $k \in I$, and let $u\colon \bk \to A$ be the unit map. Using the isomorphism $A^{\otimes(I\setminus\{k\})}\cong A^{\otimes(I\setminus\{k\})}\otimes \bk$, define
\[
\eta_k\colon A^{\otimes(I\setminus\{k\})}\longrightarrow A^{\otimes I}
\]
as the composition $A^{\otimes(I\setminus\{k\})}\cong A^{\otimes(I\setminus\{k\})}\otimes \bk \xrightarrow{\ \mathrm{id}\otimes u\ }A^{\otimes(I\setminus\{k\})}\otimes A\to A^{\otimes I}$. Set $V_k:=\operatorname{Im}(\eta_k)\subset A^{\otimes I}$. Thus, $V_k$ consists of the tensors whose $k$--leg is trivial.
\end{definition}

We now formalize a particular type of \emph{multiplication} in which two specified legs are multiplied, and the result is written in a third leg that is (formerly) the unit.

\begin{definition}[Multiplication]\label{def:contraction}
Let $i,j,k\in I$ be pairwise distinct indices. The \textit{multiplication along $(i,j)$ into $k$} is the $\bk$-linear map
\[
\mathbf{m}_{i,j \rightarrow k} \colon V_k \longrightarrow A^{\otimes I}
\]
defined by
\begin{equation}\label{eq:contraction-rule}
\mathbf{m}_{i,j\rightarrow k}\big((a)_i(b)_jx\big)\ :=\ (ab)_kx
\qquad a,b\in A,\ \ x\ \text{supported on }I\setminus\{i,j,k\},
\end{equation}
and extended linearly. In other words: for every elementary tensor whose $k$--legs is equal to the unit, the map replaces the elements in the $i$-- and $j$--legs by the unit and places the product of these elements in the $k$--leg of the output.
\end{definition}
\begin{notation}
  In what follows, we will often write composition as $f \circ g = g \sslash f$.
\end{notation}

\begin{example}\label{ex:two-contractions}
Let $I=\{i,j,k,i',j',k'\}\sqcup R$ with the six indices pairwise distinct. For $a,b,a',b'\in A$ and $x\in A^{\otimes R}$,
\[
\big(a_i b_ja'_{i'} b'_{j'}x\big) \sslash  \mathbf{m}_{i',j'\to k'} \sslash \mathbf{m}_{i,j\to k}
 =
(ab)_k(a'b')_{k'}x
 =
\mathbf{m}_{i',j'\to k'}\circ \mathbf{m}_{i,j\to k}\big(a_i b_ja'_{i'} b'_{j'}x\big).
\]
\end{example}
Lastly, we define the multi-tensor multiplication.
\begin{definition}\label{rem:variants}
For pairwise distinct $i_1,\dots,i_m,k$ and $x\in V_k$, we define the \textit{$m$--ary multiplication}
        \begin{align*}
&x \sslash \mathbf{m}_{i_1,\dots,i_m\to k}
\ :=\ \\
&\text{ replace the $i_r$--legs by 1 and write }a_{i_1}\cdots a_{i_m}\text{ into the $k$--leg.}
\end{align*}
\end{definition}
\begin{example}
  Let $I = \{1,2,3,4,5\}$ and $a,b,c,d \in A$. Then
  \[
    a_{1}b_{2}c_{3}d_{5} \sslash \mathbf{m}_{1,2,3 \to 4} = (a b c)_{4} d_{5}.
    \]
  \end{example}

\section{The Heisenberg Algebra}
\label{sec:heisenberg}
\noindent The objective of this section is to equip the reader with a range of tools to manage the occurrences of the Heisenberg algebra. This includes methods for normal ordering operators and an approach to multiplying functions of normally ordered operators. We hereby employ the theory of Gaussian integrals, in line with \cite{bar2021perturbed}.

\subsection{Normal ordering and Exponentials}
\label{sec:normord}
Let $\bk$ be a field of characteristic zero and $n \in \mathbb{N}$. Moreover, let $I$ be an ordered set with $n$ elements with ordering $i_{1} < i_{2}< \dots < i_{n}$.
\begin{definition}
  \label{def:heisenberg}
The \emph{$n$-th Weyl--Heisenberg algebra} $A_n(\bk)$ over $\bk$ is the unital associative $\bk$-algebra generated by $2n$ elements $\{\mathbf{x}_i,\mathbf{p}_i\}_{i \in I}$ subject to the relations
\begin{align*}
\mathbf{p}_i \mathbf{p}_j - \mathbf{p}_j \mathbf{p}_i = 0, \qquad \mathbf{x}_i \mathbf{x}_j - \mathbf{x}_j \mathbf{x}_i = 0, \qquad \mathbf{p}_i \mathbf{x}_j - \mathbf{x}_j \mathbf{p}_i = \delta_{ij}
\end{align*}
for all $i,j \in I$, where $\delta_{ij}$ is the Kronecker delta.
\end{definition}
\begin{remark}
The relations in the definition of the Weyl--Heisenberg algebra can also be expressed using the commutation bracket. For any two elements $\mathbf{x}$ and $\mathbf{y}$, their commutation bracket is defined as $[\mathbf{x}, \mathbf{y}] := \mathbf{x}\mathbf{y} - \mathbf{y}\mathbf{x}$. With this notation, the defining relations of the $n$-th Weyl--Heisenberg algebra can be written as
\begin{align*}
[\mathbf{p}_i, \mathbf{p}_j] = 0, \qquad [\mathbf{x}_i, \mathbf{x}_j] = 0, \qquad [\mathbf{p}_i, \mathbf{x}_j] = \delta_{ij}.
\end{align*}
\end{remark}

The Weyl--Heisenberg algebra $A_n(\bk)$ is closely related to the Heisenberg Lie algebra. In terms of the symplectic basis,  the Heisenberg Lie algebra in $n$ dimensions is defined as the Lie algebra $\mathfrak{h}_n$ with generators $\{\mathbf{p}_{k}\}_{k=1}^n \cup \{\mathbf{x}_{k}\}_{k=1}^n \cup \{\mathbf{z}\}$ and relations
\[
[\mathbf{p}_i, \mathbf{p}_j] = 0, \qquad [\mathbf{x}_i, \mathbf{x}_j] = 0, \qquad [\mathbf{p}_i, \mathbf{x}_j] = \delta_{ij}\mathbf{z}, \qquad \text{for all $1 \leq i, j \leq n$,}. 
\]
Here, $\mathbf{z}$ is a central element. We write $U(\mathfrak{h}_n)$ for the universal enveloping algebra of $\mathfrak{h}_n$, which is a noncommutative $\bk$-algebra generated by the elements of $\mathfrak{h}_n$. The $n$-th Weyl--Heisenberg algebra $A_n(\bk)$ is isomorphic to the quotient of $U(\mathfrak{h}_n)$ by the two-sided ideal generated by $\mathbf{z} - 1$, i.e., $A_n(\bk) \cong U(\mathfrak{h}_n) / (\mathbf{z} - 1)$.

The relation to the universal enveloping algebra allows us to invoke the Poincar\'e--Birkoff--Witt theorem. We hereby endow $A_n(\bk)$ with an ordering
\begin{equation}
\label{eq:ordering}
  \mathbf{x}_{i_{1}} < \dots < \mathbf{x}_{i_{n}} < \mathbf{p}_{i_{1}} < \dots < \mathbf{p}_{i_{n}}.
\end{equation}
Any ordering of this form will be called a \textit{normal ordering}. Applying the Poincaré--Birkhoff--Witt theorem to $A_n(\bk)$, we deduce that the normally ordered monomials
\[
  x^{\alpha}p^{\beta}
  := x_{i_1}^{\alpha_{1}}\cdots x_{i_n}^{\alpha_{n}}
     p_{i_1}^{\beta_{1}}\cdots p_{i_n}^{\beta_{n}},
  \qquad (\alpha,\beta) \in \mathbb{Z}_{\geq 0}^{n} \times \mathbb{Z}_{\geq 0}^{n},
\]
form a topological $\bk \llbracket h\rrbracket$-basis of $A_n(\bk)$. From this we conclude that the \emph{normal ordering map}
\[
  \mathcal{N}_{I}\colon
  \Bigl(\bigotimes_{i \in I}\bk[x_i,p_i]\Bigr)\llbracket h \rrbracket
  \xrightarrow{\ \sim\ } A_n(\bk),
  \qquad
  \prod_{i \in I} x_i^{\alpha_{i}}p_i^{\beta_{i}}
  \longmapsto
  \Bigl(\prod_{i \in I} \mathbf{x}_i^{\alpha_{i}}\Bigr)
  \Bigl(\prod_{j \in I} \mathbf{p}_j^{\beta_{j}}\Bigr),
\]
defined on monomials and extended $\bk\llbracket h \rrbracket$-linearly and
$h$-adically continuously, is an isomorphism of topological
$\bk\llbracket h \rrbracket$-modules. We omit the subscript $I$ in $\mathcal{N}_{I}$ when it is clear from the context. To ease some of the notation, set $\bk[z_I]:=\bk[x_I,p_I]:=\bigotimes_{i \in I}\bk[x_i,p_i]$ with $z_i=(x_i,p_i)$.
\begin{example}
Let $I = \{1,2\}$ and $f = \mathbf{p}_1^{2} \mathbf{x}_1\mathbf{x}_{2} \in A_{2}(k)$. Then we find
\[
f = \mathbf{p}_1^{2}  \mathbf{x}_1 \mathbf{x}_2 = 2\mathbf{x}_2\mathbf{p}_1+\mathbf{x}_1\mathbf{x}_2\mathbf{p}_1^2 = \mathcal{N}( 2 x_{2}p_{1} + x_{1} x_{2} p_{1}^{2}  ).
\]
\end{example}
The space $\bk[z_I]$ can be endowed with the structure of a commutative algebra by viewing it as a polynomial ring. In contrast, $A_n(\bk)$ is a noncommutative algebra. Therefore, the map $\mathcal{N} \colon \bk[z_I] \to A_n(\bk)$ is not an isomorphism of algebras. We will use this distinction throughout this chapter.

As we will see in the upcoming discussion, we apply exponential and logarithmic maps repeatedly to elements of the Weyl--Heisenberg algebra. Given the significance of these operations in our approach, we formally define the deformation of the underlying algebraic structure that makes them well-defined.

\begin{definition}
Let $A$ be an algebra over a field $\bk$. The \emph{trivial deformation} of $A$ into an algebra over $\bk\llbracket h \rrbracket$, is the $\bk \llbracket h \rrbracket$-algebra that arises by extending $\bk$ to $\bk \llbracket h \rrbracket$. In other words: it is the tensor product $\bk \llbracket h \rrbracket \otimes_{\bk} A$ of $\bk \llbracket h \rrbracket$ and $A$ as algebras, where the multiplication is given by $(a \otimes x)(b \otimes y) = ab \otimes xy$ for all $a, b \in \bk \llbracket h \rrbracket$ and $x, y \in A$.
\end{definition}
Every reference to the Weyl--Heisenberg algebra in our discussion means that we consider $A_n(\bk)$ as a $\bk \llbracket h \rrbracket$-algebra. In this way, the exponential map and logarithm are well-defined, and it is possible to invoke the BCH formula.

\begin{notation}
  Assume the generators $\{\mathbf{x}_i,\mathbf{p}_i\}_{i \in I}$ of $A_n(\bk)$ to be ordered as in \eqref{eq:ordering}. Moreover, set $x = (x_{i_{1}}, \dots, x_{i_{n}})$ and $p = (p_{i_{1}}, \dots, p_{i_{n}})$. Similarly, we denote $\mathbf{x}:=(\mathbf{x}_{i_{1}},\dots,\mathbf{x}_{i_{n}})$ and $\mathbf{p} := (\mathbf{p}_{i_{1}},\dots,\mathbf{p}_{i_{n}})$. Extending the map $\mathcal{N}$ entry-wise yields $\mathbf{x}= \mathcal{N}(x)$ and $\mathbf{p} = \mathcal{N}(p)$. Throughout this chapter we shall use the Einstein summation convention.
\end{notation}

\begin{lemma}
  \label{lem:expconj}
Let $A\in \mathrm{Mat}_{n}(\bk\llbracket h \rrbracket)$. The following identities hold:
\begin{enumerate}[\normalfont i)]
  \item\label{item:1}   $\exp(h \mathbf{x}^{\top} A \mathbf{p})\mathbf{x}_{j} \exp(-h \mathbf{x}^{\top} A \mathbf{p})= \mathbf{x}_{i} (e^{hA})_{ij}$;
        \item  $\exp(h \mathbf{x}^{\top} A \mathbf{p})\mathbf{p}_{j} \exp(-h \mathbf{x}^{\top} A \mathbf{p})= \mathbf{p}_{i} (e^{-hA})^{\top}_{ij}$.
\end{enumerate}

\end{lemma}
\begin{proof}
First note that
\begin{align*}
[\mathbf{x}_{\alpha} A_{\alpha \beta} \mathbf{p}_{\beta}, \mathbf{x}_{\mu}] = \mathbf{x}_{\alpha} A_{\alpha \beta} \delta_{\beta \mu} = \mathbf{x}_{\alpha} A_{\alpha \mu}.
\end{align*}
This implies that
\begin{align*}
  \mathrm{ad}_{\mathbf{x}_{\alpha} A_{\alpha \beta} \mathbf{p}_{\beta}}^n \mathbf{x}_{\mu}= \mathbf{x}_{\alpha} (A^n)_{\alpha \mu},
\end{align*}
which yields
\begin{align*}
\exp(h \mathbf{x}_{\alpha} A_{\alpha \beta} \mathbf{p}_{\beta})\mathbf{x}_{\mu} \exp(-h \mathbf{x}_{\alpha} A_{\alpha \beta} \mathbf{p}_{\beta})= \mathbf{x}_{\alpha} (e^{hA})_{\alpha \mu},
\end{align*}
proving the first identity.

For the second identity, we use the first one. Using $[\mathbf{x}_{i},\mathbf{p}_{j}] = -\delta_{ij}$, we deduce
\[
  \exp(h \mathbf{p}^{\top} A^{\top} \mathbf{x})\mathbf{x}_{j} \exp(-h \mathbf{p}^{\top} A^{\top} \mathbf{x})= \mathbf{x}_{i} (e^{hA})_{ij}.
  \]
  We then apply the automorphism
  \[
    \mathbf{x}_{i} \mapsto \mathbf{p}_{i}, \qquad \mathbf{p}_{i} \mapsto -\mathbf{x}_{i}
    \]
    to both sides of the equation. This amounts to a Fourier transform. From this, we deduce the desired equality.
\end{proof}
Let $P \colon A_{n}(\bk) \otimes A_{n}(\bk) \to A_{n}(\bk) \otimes A_{n}(\bk)$ be the interchanging map defined by $a \otimes b \mapsto b \otimes a$. Also, let $\widetilde{P}$ be the map that interchanges the columns of a $2\times 2$ matrix.
\begin{proposition}
  \label{pro:adexp}
  Let $n=2$ and $A \in \mathrm{Mat}_{2}(\bk \llbracket h \rrbracket)$. The following identity holds:
  \[
    P \circ \mathrm{Ad}_{\exp(h \mathbf{x}^{\top} A \mathbf{p})} \mathbf{x} =  \widetilde{P}((e^{hA})^{\top}) \mathbf{x}.
    \]
\end{proposition}
\begin{proof}
This is an immediate consequence of Lemma \ref{lem:expconj}.
\end{proof}
The following proposition in inspired by \cite[Lemma 10]{bar2021perturbed}
\begin{proposition}
  \label{pro:normalordering}
Let $A\in \mathrm{Mat}_{n}(\bk \llbracket h \rrbracket)$. The following identity holds:
\[
\exp(h  \mathbf{x}^{\top} A \mathbf{p}) = \mathcal{N} \left(  \exp(x^{\top} (e^{h  A}- \mathds{1}) p)\right).
\]
\end{proposition}
\begin{proof}
Utilizing the Einstein summation convention, we establish the following identity
$$\exp(h t \mathbf{x}_{\alpha} A_{\alpha \beta} \mathbf{p}_{\beta}) = \mathcal{N} \left( \exp(x_{\alpha} (e^{h  A} - \mathds{1})_{\alpha \beta} p_{\beta}) \right).$$
To prove this proposition, we work over the field $\bk \otimes \mathbb{R}$ and show that two functions, $f$ and $g$, satisfy the same ordinary differential equation (ODE) with equal initial conditions. Define
\begin{align*}
  &f \colon \mathbb{R} \to A_n(\bk), \quad t \mapsto \exp( ht \mathbf{x}^{\top} A \mathbf{p}), \\
  &g \colon \mathbb{R} \to A_n(\bk), \quad t \mapsto
\mathcal{N} \left(  \exp(x (e^{h t  A}- \mathds{1}) p)\right).
\end{align*}
 After applying Lemma \ref{lem:expconj}, we compute the derivative of $f$
\begin{align*}
  f'(t) &= h \exp(h t \mathbf{x}_{\alpha} A_{\alpha \beta} \mathbf{p}_{\beta})\mathbf{x}_{\mu} A_{\mu \nu} \mathbf{p}_{\nu} \\
  &=h \mathbf{x}_{\alpha} (e^{ht A})_{\alpha \mu}A_{\mu \nu} f(t) \mathbf{p}_{\nu} \\
  &=h \mathbf{x}_{\mu} (e^{ht A} A)_{\mu \nu} f(t) \mathbf{p}_{\nu}.
\end{align*}
Similarly, for $g$, we find
\begin{align*}
  g'(t) &= h\mathcal{N} \left(\exp(x_{\alpha} (e^{h t A}- \mathds{1})_{\alpha \beta} p_{\beta})x_{\mu} ( A e^{h t A})_{\mu \nu} p_{\nu}\right) \\
  &=h\mathbf{x}_{\mu} ( A e^{h t A})_{\mu \nu}\mathcal{N} \left(\exp(x_{\alpha} (e^{h t A}- \mathds{1})_{\alpha \beta} p_{\beta})  \right) \mathbf{p}_{\nu}\\
  &=h\mathbf{x}_{\mu} (  e^{h t A}A)_{\mu \nu} g(t) \mathbf{p}_{\nu}.
\end{align*}
Since $A e^{h t A} = e^{h t A} A$, it follows that
\[
f'(t) = h \mathbf{x}_{\mu} (e^{h t A} A)_{\mu \nu} \mathbf{p}_{\nu} f(t) \quad \text{and} \quad g'(t) = h \mathbf{x}_{\mu} (e^{h t A} A)_{\mu \nu} \mathbf{p}_{\nu} g(t).
\]
 Thus, $f$ and $g$ satisfy the same ODE with the initial condition $f(0) = g(0) = \mathds{1}$. By the uniqueness of the solutions to ODEs, we conclude that $f(t) = g(t)$ for all $t$. Setting $t = 1$ now yields the desired result.
\end{proof}
We introduce a new notation for elements of the type for which Proposition \ref{pro:normalordering} can be applied.
\begin{notation}
  \label{def:monhom}
  Denote $\mathcal{A}_{n}(\bk) := \exp(h \mathrm{Mat}_n(\bk \llbracket h \rrbracket))$ with monoid structure given by multiplication.
  Define the map $\varphi \colon \mathcal{A}_n(\bk) \to A_n(\bk)$ by
  \[
    A \mapsto \mathcal{N} \left( \exp(  x^{\top}   (A - \mathds{1})  p ) \right).
    \]
\end{notation}

The map $\varphi$ provides a dictionary between taking tensor products and adjoining block-matrices. As shall be shown below, adjoining identity factors in a tensor product corresponds, under $\varphi$, to adjoining identity summands in a direct sum.

\begin{notation}
  Fix integers $m \ge 2$ and $1 \le i < j \le m$ and let $A\in\mathcal A_2(\Bbbk)$. Define
  \[
    \varphi_{i,j}^{m}(A)
    :=
    \bigl(I^{\otimes(i-1)} \otimes \mathrm{Id} \otimes I^{\otimes(j-i-1)} \otimes \mathrm{Id} \otimes I^{\otimes(m-j)}\bigr)\bigl(\varphi(A)\bigr).
  \]
  Here $I$ denotes the identity on the underlying tensor factor.

Let $R_k = \bk \llbracket h \rrbracket$ and $R=\bigoplus_{k=1}^m R_k$. Consider the canonical injections $\iota_k \colon R_k \to R$ and projections $\pi_k\colon R \to R_k$. Moreover, let $T\in \mathrm{End}(R_i\oplus R_j)$ having entries $T_{ab}\in \mathrm{Hom}(R_b,R_a)$ for $a,b\in\{i,j\}$.  We define the block-insertion map
  \[
    \Psi^{m}_{i,j}(T)
    :=
    \sum_{k\neq i,j}\iota_k\pi_k
    +
    \sum_{a,b\in\{i,j\}}\iota_aT_{ab}\pi_b
    \in \mathrm{End}(R).
  \]
  In other words,  $\Psi^m_{i,j}(T)$ acts as the identity on $R_k$ for $k\notin\{i,j\}$, and reduces to $T$ on $R_i\oplus R_j$. Hence, $\Psi^m_{i,j}$ places a $2\times2$ block $T$ in positions $i$ and $j$.
\end{notation}  
\begin{proposition}
  \label{pro:phipsi}
Let $A \in \mathcal{A}_{2}(\bk)$. The following equality holds:
\[
\varphi^{m}_{i,j}(A) = \varphi(\Psi^{m}_{i,j}(A)).
\]
\end{proposition}
\begin{proof}
  Let $x = (x_{1},\dots,x_{m})$ and $p = (p_{1}, \dots, p_{m})$. We have
  \begin{align*}
    \varphi(\Psi^{m}_{i,j}(A)) &= \mathcal{N} \left( \exp(  x^{\top}   (\Psi_{i,j}^{m}(A) - \mathds{1})  p ) \right) \\
    \end{align*}
    Note that
    \[
      x^{\top}   (\Psi_{i,j}^{m}(A) - \mathds{1})  p  = (x_{i}, x_{j})^{\top} (A -\mathds{1}) (p_{i}, p_{j}),
      \]
     which yields the desired result.
\end{proof}
\begin{example}
Consider
\[
  A =
  \begin{pmatrix}
    T & 0 \\
    1 - T^{2} & T
  \end{pmatrix},
  \qquad
  T = e^{ht}.
\]

When we insert this $2\times 2$ block into positions $1,2$ of a bigger $3\times 3$ direct sum via $\Psi^{3}_{1,2}$, we obtain
\[
  \Psi^{3}_{1,2}(A)
  =
  \begin{pmatrix}
    T & 0 & 0 \\[2pt]
    1 - T^{2} & T & 0 \\[2pt]
    0 & 0 & 1
  \end{pmatrix}.
\]
Note that  $\Psi^{3}_{1,2}(A)$ acts as $A$ on the first two summands and as the identity on the third.
Consistent with the previous discussion, we find
\[
  \varphi^{3}_{1,2}(A) = \varphi(A) \otimes 1 = \varphi\bigl(\Psi^{3}_{1,2}(A)\bigr),
\]
which agrees with Proposition~\ref{pro:phipsi}.
\end{example}

Combining Proposition \ref{pro:normalordering} with Lemma \ref{lem:expconj} gives rise to the following.
\begin{proposition}
  \label{pro:phixandp}
  Let $A \in \mathcal{A}_n(\bk)$. The following identities hold:
\begin{enumerate}[\normalfont i)]
    \item $\varphi(A) \mathbf{x}_j \varphi(A)^{-1} = \mathbf{x}_i A_{ij}$;
    \item $\varphi(A) \mathbf{p}_j \varphi(A)^{-1} = \mathbf{p}_i (A^{-1})^{\top}_{ij}$.
  \end{enumerate}
\end{proposition}
With the required structure in place, we can now derive the following useful identity.
\begin{theorem}
  \label{cor:identity}
  Let $A,B \in \mathcal{A}_n(\bk)$. The following equality holds:
  \[
    \varphi(A) \varphi(B)= \varphi(A B),
    \]
    i.e., $\varphi$ is a monoid homomorphism.
\end{theorem}
\begin{proof}
  First we note that $[ \mathbf{x}^{\top} A \mathbf{p},\mathbf{x}^{\top} B \mathbf{p} ]= \mathbf{x}^{\top}[A,B] \mathbf{p}$. Define $e^{h \widetilde{A}} := A$ and $e^{h \widetilde{B}} := B$. This is well-defined because the logarithm exists. We find that
\begin{align*}
\varphi(A) \varphi(B)
  =&\mathcal{N} \left( \exp( x^{\top} (A-\mathds{1}) p) \right) \mathcal{N}\left(  \exp(x^{\top}(B-\mathds{1}) p)\right)\\
  &= \exp(h\mathbf{x}^{\top} \widetilde{A} \mathbf{p}) \exp(h \mathbf{x}^{\top} \widetilde{B} \mathbf{p}) \\
                                                                                      &= \exp(h \left( \mathbf{x}^{\top} \widetilde{A} \mathbf{p} + \mathbf{x}^{\top}\widetilde{B}\mathbf{p} + \frac{h}{2} [ \mathbf{x}^{\top} \widetilde{A} \mathbf{p},\mathbf{x}^{\top} \widetilde{B} \mathbf{p} ] + \dots  \right))\\
                                                                                      &= \exp( \mathbf{x}^{\top}\left(  h \widetilde{A} + h \widetilde{B} + \frac{h^{2}}{2} [  \widetilde{A} , \widetilde{B}  ] + \dots  \right) \mathbf{p})  \\
                                                                                      &= \mathcal{N} \left(  \exp(x^{\top} \left( e^{h\widetilde{A} + h\widetilde{B} + \frac{h^2}{2} [\widetilde{A},\widetilde{B}] + \dots}-\mathds{1} \right) p)\right) \\
                                                                                      &=\mathcal{N} \left(  \exp(x^{\top} \left( AB-\mathds{1} \right) p)\right) \\
  &= \varphi(AB)
\end{align*}
as claimed.
\end{proof}

\subsection{Gaussians and Contractions}
\label{sec:gausscontract}
In the previous section, we considered the multiplication $\mathbf{m} \colon A_n(\bk) \otimes A_n(\bk) \to A_{n}(\bk)$, and derived the useful identity as stated in Theorem \ref{cor:identity}. In very specific cases, this allows us to translate multiplications in $A_n(\bk) $ to multiplications in $\bk\llbracket h \rrbracket[z_I]$. Based on \cite{bar2021perturbed}, we extend this idea further: instead of performing computations in $A_n(\bk)$, we shall perform them in $\bk\llbracket h \rrbracket [z_I] $. To do so, we need to construct a dictionary between the two rings.
More precisely: to be able to carry out computations in the commutative ring $\bk \llbracket h \rrbracket [z_{I}]$, we seek an operation $m_{i,j \rightarrow k}$ such that the following diagram commutes.
\begin{equation}
  \label{eq:commdiag}
\begin{tikzcd}
 {A_n(\bk)^{\otimes \ell}} & {A_n(\bk)^{\otimes \ell}} \\
 {k\llbracket h \rrbracket [z_{I}]} & {k\llbracket h \rrbracket [z_{I}]}
 \arrow["{\mathcal{N}}", from=2-1, to=1-1]
 \arrow["{\mathbf{m}_{i,j\rightarrow k}}", from=1-1, to=1-2]
 \arrow["{m_{i,j \rightarrow k}}"', from=2-1, to=2-2]
 \arrow["{\mathcal{N}}"', from=2-2, to=1-2]
\end{tikzcd}
\end{equation}

Let $\bk$ be a field containing $\mathbb{Q}$. Fix an integer $n\ge 1$ and a finite set $I$.
Define
\[
\mathcal{B} := \bk\llbracket h\rrbracket[z_I].
\]
We introduce additional variables $r=(r_1,\dots,r_n)$ and $s=(s_1,\dots,s_n)$, by considering the $\mathcal{B}$-algebra defined as
\[
R := \mathcal{B}\llbracket r_1,\dots,r_n,s_1,\dots,s_n\rrbracket.
\]
Partial derivatives $\partial_{r_i},\partial_{s_i}$ are regarded as continuous $\mathcal{B}$-linear derivations of $R$.

\begin{definition}
Define the $\mathcal{B}$-linear map
\[
\big\langle\,\cdot\,\big\rangle_{r,s}\colon R\longrightarrow \mathcal{B}
\]
by
\begin{equation}\label{eq:contraction}
  \big\langle f \big\rangle_{r,s}
  :=
  \Big[\exp \Big(\sum_{i=1}^n \frac{\partial^2}{\partial r_i\partial s_i}\Big)f\Big]_{r=s=0}.
\end{equation}
Equivalently, if $\alpha=(\alpha_1,\dots,\alpha_n)\in\mathbb{N}^n$ and
$\partial_r^\alpha=\prod_i \partial_{r_i}^{\alpha_i}$, $\partial_s^\alpha=\prod_i \partial_{s_i}^{\alpha_i}$, then
\[
\big\langle f \big\rangle_{r,s}
=
\sum_{\alpha\in\mathbb{N}^n}\frac{1}{\alpha!}
\Big(\partial_r^\alpha\partial_s^\alpha f\Big)\Big|_{r=s=0},
\qquad
\alpha!\coloneqq \prod_{i=1}^n \alpha_i! .
\]
\end{definition}

\begin{remark}
The reader may be familiar with the language of physics, where the map $\langle \, \cdot \, \rangle_{r,s}$ is uniquely determined by
$\langle 1\rangle_{r,s}=1$, $\langle r_i s_j\rangle_{r,s}=\delta_{ij}$, $\langle r_i r_j\rangle_{r,s}=0$,
$\langle s_i s_j\rangle_{r,s}=0$, and Wick’s theorem (pairwise contraction).
It corresponds to performing a “Gaussian integral”
\[
\langle g(r,s)\rangle_{r,s}\sim\int e^{-\sum_i r_i s_i}g(r,s){\mathrm d}r{\mathrm d}s,
\]
but we work purely formally via \eqref{eq:contraction}, avoiding any convergence issues.
\end{remark}

For a real, symmetric, positive-definite matrix $A$, the standard Gaussian integral is
\[
  \int_{\mathbb{R}^{n}} \exp \left(-\tfrac{1}{2} x^{\top} A x\right) \mathrm{d}x
  = (2\pi)^{n/2} (\det A)^{-1/2}.
\]
When polynomial perturbations are inserted, the resulting higher-order perturbations can be computed using Wick's theorem. In the framework of the contractions introduced above, the following analogous identity holds.
\begin{theorem}[Contraction Theorem]\label{thm:contraction}
Let $f,g\in \mathcal{B}^{n}$ and $W \in \bk \llbracket h \rrbracket^{n \times n}$. Write
\[
gs\coloneqq \sum_{i=1}^n g_i s_i,\qquad
rf\coloneqq \sum_{i=1}^n r_i f_i,\qquad
rWs\coloneqq \sum_{i,j=1}^n r_i W_{ij} s_j .
\]
Assume that the matrix $\mathbf 1 - W$ is invertible.
Set $\widetilde W\coloneqq (\mathbf 1-W)^{-1}$, and let $P$ be a polynomial in $(r,s)$. One has the identity
\begin{align*}
&\left\langle P(r,s)\exp \big(gs + rf + rWs\big) \right\rangle_{r,s} \\
&=
\det(\widetilde W)\exp \big(g\widetilde Wf\big)
\left\langle P\big(r+g\widetilde W,\ \widetilde W(s+f)\big) \right\rangle_{r,s}.
  \end{align*}
Here, $g\widetilde W f \in \mathcal{B}$ denotes the scalar $\sum_{i,j} g_i\widetilde W_{ij}f_j$.
\end{theorem}
\begin{proof}
  See \cite[Section 2]{bar2021perturbed}.
\end{proof}

The notion of contraction is particularly useful in an algebraic setting.
The following discussion is fundamentally based on the theory of generating functions.

\begin{theorem}
  \label{thm:multiplication}
  Let $r = (x_{i},x_{j}, \pi_{i}, \pi_{j})$ and $s = (\xi_{i},\xi_{j},p_{i},p_{j})$.
  The map
  \begin{align*}
    &m_{i,j \rightarrow k} \colon k\llbracket h \rrbracket [z_{I}] \to k\llbracket h \rrbracket [z_{I}], \\
  &f \mapsto \left\langle  e^{(\pi_{i} + \pi_{j})p_{k} + (\xi_{i} + \xi_{j})x_{k} + \pi_{i}\xi_{j}}f \right\rangle_{r,s}
    \end{align*}
renders diagram \eqref{eq:commdiag} commutative.
\end{theorem}
\begin{proof}
  See \cite[Section 2]{bar2021perturbed}.
\end{proof}

\begin{example}
  Let $\mathbf{f} = \mathbf{x}_{i}\mathbf{p}_{j}$. Clearly,
  \[
    \mathbf{m}_{i,j \rightarrow k}(\mathbf{f}) = \mathbf{x}_{k} \mathbf{p}_{k}.
  \]

  Next, consider $\mathcal{N}^{-1}(\mathbf{f}) = x_{i}p_{j} =: f$.
  Define $r =  (x_{i},x_{j}, \pi_{i}, \pi_{j})$ and $s = (\xi_{i},\xi_{j},p_{i},p_{j})$. We compute
  \begin{align*}
    \left\langle e^{(\pi_{i} + \pi_{j})p_{k} + (\xi_{i} + \xi_{j})x_{k} + \pi_{i}\xi_{j}} f \right\rangle_{r,s}
    &= \left. e^{(\partial_{p_{i}} + \partial_{p_{j}})p_{k} + (\partial_{x_{i}} + \partial_{x_{j}})x_{k} + \partial_{p_{i}} \partial_{x_{j}}} x_{i}p_{j} \right|_{s=0} \\
    &= x_{k}p_{k}.
  \end{align*}

  Now, let $\mathbf{g} = \mathbf{p}_{i}\mathbf{x}_{j}$. In this case,
  \[
    \mathbf{m}_{i,j \rightarrow k}(\mathbf{g}) = \mathbf{p}_{k}\mathbf{x}_{k}
    = \mathbf{x}_{k}\mathbf{p}_{k} + 1.
  \]
  Since $\mathcal{N}^{-1}(\mathbf{g}) = p_{i}x_{j} =: g$, we obtain
  \begin{align*}
    \left\langle e^{(\pi_{i} + \pi_{j})p_{k} + (\xi_{i} + \xi_{j})x_{k} + \pi_{i}\xi_{j}} g \right\rangle_{r,s}
    &= \left. e^{(\partial_{p_{i}} + \partial_{p_{j}})p_{k} + (\partial_{x_{i}} + \partial_{x_{j}})x_{k} + \partial_{p_{i}} \partial_{x_{j}}} p_{i}x_{j} \right|_{s=0} \\
    &= x_{k}p_{k} + 1.
  \end{align*}
 These computations are in agreement with Theorem~\ref{thm:multiplication}.
\end{example}

\begin{example}\label{ex:stitching}
Set $r=(x_{1},x_{2},\pi_{1},\pi_{2})$ and $s=(\xi_{1},\xi_{2},p_{1},p_{2})$. Then
\begin{align*}
\varphi \left(  \begin{pmatrix}
\alpha & \beta & \theta \\
\gamma & \delta & \varepsilon \\
\phi & \psi & \Xi
\end{pmatrix} \right)
\sslash \mathbf{m}_{1,2 \rightarrow \bar{2}}
=
\mathcal{N} \left(\langle \exp(gs + rf + r W s)\rangle_{r,s}\right),
\end{align*}
with
\[
W=
\begin{pmatrix}
0 & 0 & \alpha-1 & \beta \\
0 & 0 & \gamma & \delta-1 \\
0 & 1 & 0 & 0 \\
0 & 0 & 0 & 0
\end{pmatrix},
\]
and
\[g=(x_{\bar{1}},x_{\bar{1}},\phi x_{3},\psi x_{3}),
f=(\theta p_{3},\varepsilon p_{3},p_{\bar{1}},p_{\bar{1}}).\]
Let $\widetilde{W}\coloneqq(\mathbb{1}-W)^{-1}$. Using Theorem~\ref{thm:contraction}, we obtain
\begin{align*}
\varphi \left(  \begin{pmatrix}
\alpha & \beta & \theta \\
\gamma & \delta & \varepsilon \\
\phi & \psi & \Xi
\end{pmatrix} \right)
\sslash \mathbf{m}_{1,2 \rightarrow \bar{2}}
& \overset{\bar{2} \mapsto 2}{=}
\frac{1}{1-\gamma} 1 \otimes
\varphi \left(  \begin{pmatrix}
\displaystyle \beta + \frac{\alpha \delta}{1-\gamma} &
\displaystyle \theta + \frac{\alpha \varepsilon}{1-\gamma} \\
\displaystyle \psi + \frac{\delta \phi}{1-\gamma} &
\displaystyle \Xi + \frac{\phi \varepsilon}{1-\gamma}
\end{pmatrix} \right),
\end{align*}
where in the last step we relabel $\bar{2}\mapsto 2$.
\end{example}

\begin{remark}
The computation in Example~\ref{ex:stitching} connects neatly to the notation used in Gassner calculus; see~\cite{vo2017metamonoidsfoxmilnorcondition} for details.
\end{remark}
 \begin{remark}
Note that in Example~\ref{ex:stitching}, a contraction produces a scalar factor $\frac{1}{1-\gamma}$. This is a general feature of contractions, as is clear from Theorem~\ref{thm:contraction}. After performing an arbitrary number of contractions on $\varphi(P)$ with $P \in \mathcal{A}_{n}(\bk)$, $n \in \mathbb{N}$, we end up with
\[
  Z := \omega \cdot \varphi(Q) \qquad \omega \in \bk \llbracket h \rrbracket, \quad Q \in \mathcal{A}_{m}(\bk), \text{ with } m \in \mathbb{N}.
\]   We refer to the value $\omega$ as the \emph{scalar part} of $Z$, while $Q$ is referred to as the \emph{matrix part} of $Z$.
  \end{remark}

Although Theorem \ref{thm:multiplication} allows us to perform multiplications between two tensor factors, it is desirable to be able to compute multiplications in bulk. The following two propositions allow us to do exactly that.

\begin{proposition}
  \label{thm:gencont}
  Let $n\in \mathbb{N}$ and let 
$r = \left(x, \pi\right) \in R,$
and
$s = \left(\xi, p\right) \in R$ with
\begin{align*}
  x=(x_{1},\dots,x_{n}), \quad \pi = (\pi_{1},\dots \pi_{n}) \quad  s = (\xi_{1},\dots, \xi_{n}) \quad p = (p_{1},\dots,p_{n}).
\end{align*}
Also, define
\[
  g = \sum_{i=1}^{n} x_{\bar n} e_i \in \mathcal{B}^{2n}[x_{\bar{n}},p_{\bar{n}}],
  \qquad
  f = \sum_{i=n+1}^{2n} p_{\bar n} e_i \in  \mathcal{B}^{2n}[x_{\bar{n}},p_{\bar{n}}],
\]
where $e_i$ are the standard basis vectors.
Define the matrix $ M_{n} $ as follows:
$$M_{n} = \sum_{1 \leq i < j \leq n} E_{ij},$$
where $ E_{ij} $ is a matrix with $1$ in the $ (i,j) $ position and $0$ elsewhere. Consider the map
$$m_{1,2,\ldots,n \to \bar{n}} \colon \bk \llbracket h \rrbracket[z_{I}] \to \bk \llbracket h \rrbracket[z_{I}]$$
defined by
$$ H \mapsto  \left\langle e^{\pi M_{n} \xi + rf +gs} H \right\rangle_{r,s}.$$ The following equality holds:
$$\mathbf{m}_{1,2,\ldots,n \to \bar{n}} \circ \mathcal{N} = \mathcal{N} \circ m_{1,2,\ldots,n \to \bar{n}}.$$
\end{proposition}
\begin{proof}
  See \cite[Section 2 \& 3]{bar2021perturbed}.
\end{proof}
\begin{proposition}[Contracting in bulk]
  \label{pro:contrbulk}
Let
$A \in \mathcal{A}_{n}(\bk)$. Denote by $\hat{A} \in \mathcal{A}_{n-1}(\bk)$
the matrix obtained from $A$ by deleting its first row and its last column.
Then the scalar part of $  \varphi(A) \sslash \mathbf{m}_{1,2,\dots,n \to \bar{n}}$ is equal to
\begin{align*}
  \frac{1}{\Det(\mathds{1} - \hat{A})}.
  \end{align*}
\end{proposition}
\begin{proof}
Let $M_{n}$ be as in Theorem~\ref{thm:gencont}, and define
\begin{equation}
    \label{eq:doubleu}
  W =
  \begin{pmatrix}
    0 & A - \mathds{1} \\
    M_{n} & 0
  \end{pmatrix}.
\end{equation}
To compute
\[
  \varphi(A) \sslash \mathbf{m}_{1,2,\dots,n \to \bar{n}}
\]
we apply Proposition~\ref{thm:gencont} together with the contraction theorem.
Observe that
\[
  \Det\bigl(\mathds{1} - W\bigr)
  = \Det
    \begin{pmatrix}
      \mathbb{I}_{n} &  \mathbb{I}_{n}-A \\
       - M_{n} & \mathbb{I}_{n}
    \end{pmatrix}.
\]
Define the $(n-1) \times (n-1)$-matrix
    \[
N =
\begin{pmatrix}
0 & 1 & 0 & \cdots & 0 \\
0 & 0 & 1 & \ddots & \vdots \\
0 & 0 & 0 & \ddots & 0 \\
\vdots & \vdots & \vdots & \ddots & 1 \\
0 & 0 & 0 & \cdots & 0
\end{pmatrix}.
\]

We expand the determinant $\Det (1-W)$ along the last row, which contains zeros everywhere except in the last entry, so only one cofactor term remains. We do the same for the first column, which contains zeros everywhere except in the first entry.  Using the fact that $N + M_{n-1}N = M_{n-1}$, we find
\begin{align*}
\Det \left(\mathds{1} - W \right) &= \Det \left( \mathbb{I}_{2n-2} - \begin{pmatrix} 0 & \hat{A} - N \\ \mathbb{I}_{n-1} + M_{n-1} & 0 \end{pmatrix} \right) \\
  &=\Det \left( \mathbb{I}_{n-1} - ( \mathbb{I}_{n-1} + M_{n-1})(\hat{A} - N) \right) \\
  &= \Det \left( \mathbb{I}_{n-1} - \hat{A} -M_{n-1}\hat{A} + N + M_{n-1}N \right)\\
  &= \Det \left( \mathbb{I}_{n-1} - \hat{A} +  M_{n-1}(\mathbb{I}_{n-1} - \hat{A}) \right) \\
                                  &= \Det ((M_{n-1} + \mathbb{I}_{n-1})( \mathbb{I}_{n-1}-\hat{A})) \\
  &= \Det(\mathbb{I}_{n-1} - \hat{A}  ),
 \end{align*}
 where in the last step we used $\Det(M_{n-1}+\mathbb{I}_{n-1}) = 1$.
 \end{proof}

\section{The Alexander Polynomial}
\label{sec:alexander}

\noindent In this section, we construct a knot invariant using the Heisenberg algebra $\A_{1}(\bk)$. We use the universal invariant, XC-algebras, and twisting; all of which we briefly recall below. For a formal treatment, see \cite[Secs.~2.2–2.4]{bosch2025largecolorexpansionderiveduniversal}.

A standard approach for constructing knot invariants utilizes ribbon Hopf algebras~\cite{ReshetikhinTuraev1990,TuraevQI}. Informally, a \emph{ribbon Hopf algebra} is a Hopf algebra $A$ equipped with
distinguished invertible elements
\[
R=\sum_i \alpha_i\otimes \beta_i \in A\otimes A,\qquad \kappa\in A^\times,
\]
which turns the tensor category $\mathrm{fMod}_A$ of finite-dimensional $A$-modules into a ribbon  category. Here $\kappa$ denotes the balancing element.

It is equally possible to construct knot invariants universally. This is thoroughly described in \cite[Section 4.1]{ohtsuki2002quantum}. We recall the essentials, using a slightly different convention.

Let $(A,R,\kappa)$ be a ribbon Hopf algebra. Consider a long, framed, oriented knot $\mathcal{K}$ expressed within a diagram in which every strand of every crossing is locally oriented upward. Denote $R=\sum_i \alpha_i\otimes \beta_i \in A\otimes A$ and $R^{-1} = \sum \overline{\alpha_i}\otimes \overline{\beta}_i \in A\otimes A$. To each generator of the knot diagram, we place beads in accordance with the following rules:
\begin{equation*}
\centre{
\centering
\includegraphics[width=0.8\textwidth]{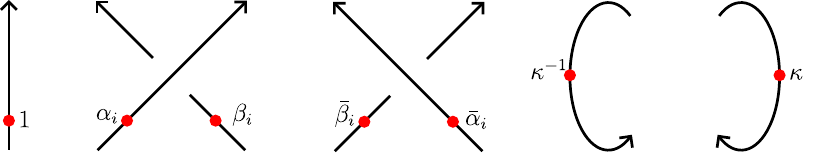}}
\end{equation*}

Traversing the knot once in the given orientation and multiplying the bead labels from right to left produces an element
\[
Z_A(\mathcal{K})\in A,
\]
the \textit{universal invariant} associated to $\mathcal{K}$.

\begin{example}
  Let $\mathcal{K}$ be the trefoil as illustrated below. We have already decorated the crossings and closures as described above.
\begin{figure}[H]
\centering
\begin{tikzpicture}[scale=1]
  \node[anchor=south west, inner sep=0] (img) at (0,0)
    {\includegraphics[scale=0.7]{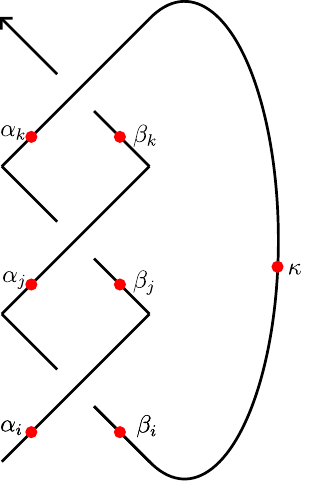}};

\end{tikzpicture}
\end{figure}
The universal invariant, associated with the trefoil, is given by
\[
  Z(\mathcal{K})= \sum_{i,j,k}  \beta_k\alpha_j \beta_i\kappa \alpha_k\beta_{j }\alpha_i.
  \]
\end{example}
To obtain a knot invariant, one does not need a full ribbon Hopf algebra structure on an algebra $A$. Instead, Becerra showed in \cite{becerra2025refinedfunctorialuniversaltangle} that it suffices to consider an algebra structure on $A$ along with specific relations between the elements $R$ and $\kappa$ that ensure the invariance of $Z_A(\mathcal{K})$ under changes of the knot diagram that do not alter the isotopy type of the knot $\mathcal{K}$ it represents. Equivalently, elements $R$ and $\kappa$ satisfy certain relations that force invariance of $Z_A(\mathcal{K})$ under the framed Reidemeister moves. This more general structure, given by the triple $(A, R, \kappa)$, is referred to as an \textit{XC-algebra}. Formally, it is defined as follows.
\begin{definition}
Let $ \bk $ be a commutative ring with a unit, and let $(A, \mu, 1) $ be a $ \bk $-algebra. An \emph{$XC $-structure} on $ A $ consists of two invertible elements,
\[
    R \in A \otimes A, \qquad \kappa \in A,
\]
called the \emph{universal $R $-matrix} and the \emph{balancing element}, respectively, which must satisfy the following conditions:
\begin{enumerate}
    \item\label{item:1} $ R = (\kappa \otimes \kappa) \cdot R \cdot (\kappa^{-1} \otimes \kappa^{-1}) $,
  \item $\mu^{[3]}(R_{13} \cdot \kappa_2) = \mu^{[3]}(R_{31} \cdot \kappa_2^{-1}) $,
        \item $1 \otimes \kappa^{-1} = (\mu \otimes \mu^{[3]})(R_{15} \cdot R_{23}^{-1} \cdot \kappa_{4}^{-1}),$
    \item $ \kappa \otimes 1 = (\mu^{[3]} \otimes \mu)(R_{34}^{-1} \cdot R_{15} \cdot \kappa_2) $,
    \item $R_{12} R_{13} R_{23} = R_{23} R_{13} R_{12} $,
\end{enumerate}
where $ \mu^{[3]} $ denotes the three-fold multiplication map. For indices $1 \leq i, j \leq n $ with $ i \neq j $, we define
\[
    R_{ij} := \begin{cases}
        (1^{\otimes (i-1)} \otimes \mathrm{Id} \otimes 1^{\otimes (j-i-1)} \otimes \mathrm{Id} \otimes 1^{\otimes (n-j)})(R^{\pm 1}), & i > j, \\
        (1^{\otimes (j-1)} \otimes \mathrm{Id} \otimes 1^{\otimes (i-j-1)} \otimes \mathrm{Id} \otimes 1^{\otimes (n-i)})(\mathrm{flip}_{A,A}R^{\pm 1}), & j > i.
    \end{cases}
\]
Similarly, we write $\kappa_i^{\pm 1} = (1^{\otimes (i-1)} \otimes \mathrm{Id} \otimes 1^{\otimes (n-i)})(\kappa^{\pm 1}) $.

A triple $ (A, R, \kappa) $ consisting of a $\bk $-algebra equipped with an $ XC $-structure is called an \emph{$XC $-algebra}.
\end{definition}

Having established the contraction formalism in the previous section, we now define the Alexander polynomial in terms of Gaussian contractions. 
\begin{propdef}
  \label{pro:alexxc}
  Let $T = e^{ht}$. The triple $\mathscr{A} = (A_{1}(\bk),R = \varphi(M),T^{-1})$ with
  \[
  M = \begin{pmatrix} T & 0 \\ 1-T^{2} & T \end{pmatrix}
\]
is an XC-algebra.
\end{propdef}
\begin{proof}
  Since the balancing element in the triple $\mathcal{A}$ is a scalar, it is obvious that the first property of XC-algebras holds. For the second property, we use Theorem \ref{thm:contraction} and \ref{thm:multiplication}. We find
  \begin{align*}
   \mu^{[3]}(R_{13} \cdot \kappa_2^{-1}) =  \varphi (M) \sslash \mathbf{m}_{1,3 \to \bar{1}} = \mathrm{det} \begin{pmatrix} 0 & M-I \\ C & 0 \end{pmatrix} = T^{-1},
\end{align*}
with
\[
  C = \begin{pmatrix} 0 & 1 \\ 0 & 0 \end{pmatrix}.
\]
Similarly, we find
  \begin{align*}
    \mu^{[3]}(R_{31} \cdot \kappa_2) =  \varphi (M) \sslash \mathbf{m}_{3,1 \to \bar{1}} = \mathrm{det} \begin{pmatrix} 0 & M^{\top}-I \\ C & 0 \end{pmatrix} = T^{-1}.
    \end{align*}
    Hence, the second property of XC-algebras holds. By an analogous computation, the third and fourth properties can be shown.

  For the last equality, we make use of Proposition \ref{pro:phipsi}. Since $R_{12}, R_{13},R_{23} \in \mathcal{A}_{1}(\bk)^{\otimes 3}$, we write
\begin{align*}
  R_{12} &= \varphi  \left(  \begin{pmatrix} T & 0 & 0 \\ 1-T^{2} & T & 0 \\ 0 & 0 & 1\end{pmatrix} \right) \quad
  R_{13} = \varphi  \left(  \begin{pmatrix} T & 0 & 0 \\ 0 & 1 & 0 \\ 1-T^{2} & 0 & T \end{pmatrix} \right) \\
  R_{23} &= \varphi \left( \begin{pmatrix} 1 & 0 & 0 \\ 0 & T & 0 \\ 0 & 1-T^{2} & T \end{pmatrix} \right).
\end{align*}
Using Theorem \ref{cor:identity}, we compute the product using matrix multiplication:
\[
  R_{12} R_{13} R_{23} = R_{23} R_{13}R_{12} = \varphi \left(  \begin{pmatrix} T^{2} & 0 & 0 \\ T-T^{3} & T^{2} & 0 \\ 1-T^{2} & T-T^{3} & T^{2} \end{pmatrix} \right). \qedhere
  \]

\end{proof}

\begin{example}
  \label{ex:firstalexander}
Let
\begin{align*}
  \mathbf{R}_{i,j} &= \mathcal{N} \left(\mathrm{Exp} \big((T-1)x_{i}p_{i} + (1-T^{2}) x_{j}p_{i} + (T-1) x_{j}p_{j}\big)\right) \\
  &= \mathcal{N} \left(\mathrm{Exp} \big(x^{\top}(M-\mathbb{I})p\big)\right)
\end{align*}
with
\[
  x = (x_{i}, x_{j}), \quad p = (p_{i},p_{j}), \quad
  M = \begin{pmatrix} T & 0 \\ 1-T^{2} & T \end{pmatrix}.
\]
Recall the interchange map
\[
  P \colon A_{n}(\Bbbk) \otimes A_{n}(\Bbbk) \to A_{n}(\Bbbk) \otimes A_{n}(\Bbbk),
  \quad a \otimes b \mapsto b \otimes a.
\]
We shall compute the universal invariant associated to the XC-algebra $\mathcal{A}$ for the trefoil. Following the prescription of the universal invariant, and applying the stitching procedure, we compute
\begin{align*}
  \mathbf{Z} := (T^{-1}\mathbf{R}_{1,4}\mathbf{R}_{5,2}\mathbf{R}_{3,6}) \sslash \mathbf{m}_{1,2 \rightarrow 1} \sslash \mathbf{m}_{4,5 \rightarrow 2} \sslash \mathbf{m}_{1,3 \rightarrow 1} \sslash \mathbf{m}_{2,6 \rightarrow 2} \sslash \mathbf{m}_{1,2 \rightarrow 1}.
\end{align*}
Instead of invoking the contraction theorem, we make use of Theorem~\ref{cor:identity}. A direct computation gives
\begin{align*}
  &(\mathbf{R}_{1,4}\mathbf{R}_{5,2}) \sslash \mathbf{m}_{1,2 \rightarrow 1} \sslash \mathbf{m}_{4,5 \rightarrow 2} = \mathbf{R}_{1,2} P(\mathbf{R}_{1,2}) \\
  &= \mathcal{N}\left(\mathrm{Exp}\big((T^{2}-1)x_{1}p_{1} + (T-T^{3})x_{2}p_{1} + (T-T^{3})x_{1}p_{2} + (T^{4}-T^{2})x_{2}p_{2}\big)\right).
\end{align*}
With this, we obtain
\begin{align*}
  &(\mathbf{R}_{1,4}\mathbf{R}_{5,2}\mathbf{R}_{3,6}) \sslash \mathbf{m}_{1,2 \rightarrow 1} \sslash \mathbf{m}_{4,5 \rightarrow 2} \sslash \mathbf{m}_{1,3 \rightarrow 1}  \sslash \mathbf{m}_{2,6 \rightarrow 2} \\
  &= \mathbf{R}_{1,2}P \big(\mathbf{R}_{1,2} P(\mathbf{R}_{1,2})\big) \\
  &= \mathcal{N} \left(\mathrm{Exp} \big(x^{\top}(B-\mathbb{I})p\big)\right)
\end{align*}
where
\[
  B = \begin{pmatrix} T - T^{3} + T^{5} & T^{2} - T^{4} \\[6pt] 1 - T^{2} + T^{4} - T^{6} & T - T^{3} + T^{5} \end{pmatrix}.
\]
Finally, let
\[
  C = \begin{pmatrix} 0 & 1 \\ 0 & 0 \end{pmatrix}.
\]
Applying Theorem~\ref{thm:contraction}, we arrive at
\begin{align*}
  \mathbf{Z} &\doteq \det \begin{pmatrix} \mathbb{I} & B - \mathbb{I} \\ C & \mathbb{I} \end{pmatrix}^{-1} \\
  &\doteq \frac{1}{T^{2} - T^{4} + T^{6}}.
\end{align*}
where $\doteq$ denotes equality up to multiplication by $\pm T^{m}$ for some integer $m$. After substituting $T \mapsto T^{1/2}$, this expression agrees with the Alexander polynomial up to rescaling.
\end{example}

\begin{figure}[H]
    \centering
    \includegraphics[scale=1]{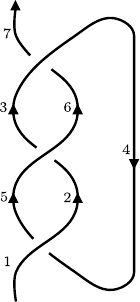}
    \caption{Labeled trefoil}
  \label{fig:labtrefoil}
\end{figure}

In the next paragraphs, we recall the definition of the Fox matrix as provided in \cite{barnatan2024perturbedalexanderinvariant}. Subsequently, we show that the universal invariant related to the XC-algebra from Proposition~\ref{pro:alexxc} gives rise to the Alexander polynomial.

Let $\mathcal{K}$ be an oriented knot with $n$ crossings. Draw it in the plane as a long knot diagram $D$ in such a way that at every crossing both strands are oriented upward (this is always possible, since we may rotate crossings as needed), and so that near its beginning and its end the knot is oriented upward as well. We call such a diagram an \emph{upward knot diagram}. See also \cite{becerra2024bar, bosch2025largecolorexpansionderiveduniversal} for a more formal definition. An example of an upward knot diagram is shown in Figure~\ref{fig:labtrefoil}.

 We now label each edge of the diagram with an integer: a running index $k$ that ranges from $1$ to $2n+1$ along the orientation of the knot. In the above example, the indices range from $1$ to $7$.

 Define the $(2n+1)\times (2n+1)$-matrix
 \begin{equation}
   \label{eq:alexandermatrix}
 A = I + \sum_{c} A_{c},
 \end{equation}
 where for each crossing $c$, the matrix $A_{c}(T)$ is zero except in the following blocks:

\begin{figure}[H]
\centering
\begin{tikzpicture}[baseline=(img.base)]

  \node[inner sep=0] (img) at (0,0) {
    \includegraphics[scale=0.8]{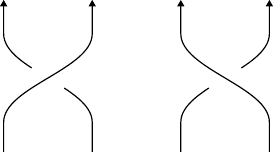}
  };

  \node at (-2, -0.8) {$i$};
  \node at (-0.8, -0.8) {$j$};
  \node at (-2.1, 0.9) {$j^{+}$};
  \node at (-0.9, 0.9) {$i^{+}$};
  \node at (0.4, -0.8) {$j$};
  \node at (1.6, -0.8) {$i$};
  \node at (0.3, 0.9) {$i^{+}$};
  \node at (1.5, 0.9) {$j^{+}$};
  \node at (-1.25, -1.5) {$s = +1$};
  \node at (1.25, -1.5) {$s = -1$};
  \node[right=1.2cm of img] (arrow) {\(\longrightarrow\)};

  \node[right=1.2cm of arrow, baseline] {
  \(
  \begin{array}{c|cc}
  A_c & \text{column } i^{+} & \text{column } j^{+} \\ \hline
  \text{row } i & -T^s & T^s - 1 \\
  \text{row } j & 0 & -1
  \end{array}
  \)
  };
\end{tikzpicture}
\end{figure}
The matrix $A$ arises from Fox calculus applied to the Wirtinger presentation of the fundamental group of the complement of $\mathcal{K}$. The \textit{Alexander polynomial} of $\mathcal{K}$ is then defined, up to multiplication by a unit in $\mathbb{Z}[T^{\pm 1}]$, as the determinant of $A$:

 \[
 \Delta_{\mathcal{K}}(T) \doteq \Det A,
 \]
where $\doteq$ denotes equality up to multiplication by $\pm T^{m}$ for some integer $m$.

Given an XC-algebra $\mathscr{A}' = (A,R,\kappa)$, it is possible to define another XC-algebra $\mathscr{B}' = (A,\check{R}, \kappa)$ via a process called twisting. Indeed, given an element $\psi$ such that $[\psi, \kappa] = 0$ and $[\psi \otimes \psi,R] = 0$, we define $\check{R} = (1 \otimes \psi^{-1}) R (\psi \otimes 1)$. As shown in \cite[Corollary 2.13]{bosch2025largecolorexpansionderiveduniversal} we then have $Z_{\mathscr{A}}(\mathcal{K}) = Z_{\mathscr{B}}(\mathcal{K})$. With this, we prove the following.

\begin{lemma}
  \label{lem:XCB}
The XC-algebra $\mathscr{A} = (A_{1}(\bk),R = \varphi(M),T^{-1})$ with
  \[
  M = \begin{pmatrix} T & 0 \\ 1-T^{2} & T \end{pmatrix}
\]
is related to the XC-algebra $\mathscr{B} = (A_{1}(\bk),\check{R} = \varphi(N),T^{-1})$ with

\begin{equation}
  \label{eq:N}
N = \begin{pmatrix} T^{2} & 0 \\ 1-T^{2} & 1 \end{pmatrix}
\end{equation}
via twisting with the element
 \[
   \psi := e^{h t x p } \in A_{1}(\bk).
 \]
\end{lemma}
\begin{proof}
We have
 \[
\psi \otimes 1 = \varphi \left(  \begin{pmatrix} T & 0 \\ 0 & 1 \end{pmatrix} \right), \quad 1 \otimes \psi = \varphi \left(  \begin{pmatrix} 1 & 0 \\ 0 & T \end{pmatrix}\right) \quad \text{and} \quad  \psi \otimes \psi = \varphi \left(  \begin{pmatrix} T & 0 \\ 0 & T \end{pmatrix} \right).
   \]
   Note that the matrix $\begin{pmatrix} T & 0 \\ 0 & T \end{pmatrix}$ commutes with $\begin{pmatrix} T & 0 \\ 1-T^{2} & T \end{pmatrix}$. Hence, $\psi$ is a twisting element. Using PropositionUsing \cite[Prop. 2.9]{bosch2025largecolorexpansionderiveduniversal}, the triple $\mathcal{B} = (A_{1}(\bk), \check{R} ,1)$ with
  \[
\check{R} = (1 \otimes \psi^{-1}) R (\psi \otimes 1)
\]
defines an XC-algebra. Note that $\check{R} = \varphi(N)$ with
\[
  N = \begin{pmatrix} T^{2} & 0 \\ 1-T^{2} & 1 \end{pmatrix}. \qedhere
\]
\end{proof}

With the XC-algebra $\mathcal{B}$ in place, we conclude this section with the following theorem.

\begin{theorem}
\label{thm:alex}
Let $\mathcal{K}$ be a long knot and let $\mathscr{A}$ be the XC-algebra defined in Proposition~\ref{pro:alexxc}. Denote by $\widetilde{\mathcal{K}}$ the closed knot obtained from $\mathcal{K}$ by trivially closing its open component, and let $\omega$ be the scalar part of $Z_{\mathscr{A}}(\widetilde{\mathcal{K}})$. Then
\[
  \omega \doteq \Delta_{\widetilde{\mathcal{K}}}(T^{2}),
\]
where $\Delta_{\widetilde{\mathcal{K}}}(T)$ is the Alexander polynomial of $\widetilde{\mathcal{K}}$ and $\doteq$ denotes equality up to multiplication by a unit $\pm T^{n}$ with $n\in\mathbb{Z}$.
\end{theorem}

\begin{proof}
  Let $m$ be the number of edges of $\mathcal{K}$. Label the edges of the long knot $\mathcal{K}$ in the opposite order as described above (with Figure \ref{fig:labtrefoil} as an example). Add one extra strand above the first strand and label it by 0. The diagram $D$ now has $m+1$ edges. Stitching all the edges together gives a diagram of the knot $\mathcal{K}$.

  With Lemma \ref{lem:XCB} and \cite[Corollary 2.13]{bosch2025largecolorexpansionderiveduniversal}, we know $Z_{\mathscr{A}}(\mathcal{K}) = Z_{\mathscr{B}}(\mathcal{K})$. Denote $\prod_c$ for the product over all the crossings in $\mathcal{K}$. The universal invariant is then of the following form
 \begin{align*}
   Z_{\mathscr{A}}(\mathcal{K}) &= Z_{\mathscr{B}}(\mathcal{K}) \\
   &\doteq \prod_{c} \varphi_{i,j}^{m+1}(N) \sslash \mathbf{m}_{0,1,\dots,m} \\
   &\doteq \varphi(\widetilde{N}) \sslash \mathbf{m}_{0,1,\dots,m},
   \end{align*}
    where $N$ is defined in Eq. \eqref{eq:N} and $\widetilde{N} \in \mathcal{A}_{m+1}(\bk)$, which, consistent with Proposition \ref{pro:phipsi}, is defined as $\widetilde{N} = E_{1,1} + E_{m+1,m+1} + \sum_{c} \widetilde{N}_{c}$ with $\widetilde{N}_c$ being zero except in the following blocks:
\begin{figure}[H]
\centering
\begin{tikzpicture}[baseline=(img.base)]

  \node[inner sep=0] (img) at (0,0) {
    \includegraphics[scale=0.8]{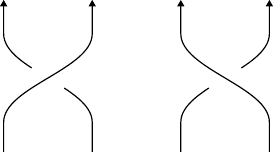}
  };

  \node at (-2, -0.8) {$i$};
  \node at (-0.8, -0.8) {$j$};
  \node at (-2.1, 0.9) {$j^{-}$};
  \node at (-0.9, 0.9) {$i^{-}$};
  \node at (0.4, -0.8) {$j$};
  \node at (1.6, -0.8) {$i$};
  \node at (0.3, 0.9) {$i^{-}$};
  \node at (1.5, 0.9) {$j^{-}$};
  \node at (-1.25, -1.5) {$s = +1$};
  \node at (1.25, -1.5) {$s = -1$};

  \node[right=0.45cm of img] (arrow) {\(\longrightarrow\)};

  \node[right=0.45cm of arrow, baseline] {
  \(
  \begin{array}{c|cc}
  \widetilde{N}_c & \text{column } i^- +1 & \text{column } j^{-}+1 \\ \hline
  \text{row } i^- +1 & T^{2s} &  0 \\
  \text{row } i^{-} +1 & 1-T^{2s} & 1
  \end{array}
  \)
  };
\end{tikzpicture}
\end{figure}
We then apply Proposition \ref{pro:contrbulk} to stitch all the edges together. It follows that the scalar part of $Z_{\mathcal{A}}(\mathcal{K})$ is equal to $\Det(\mathds{1}+ \widetilde{A})^{-1}$, where $\widetilde{A}$ is equal to the matrix $\widetilde{N}$ with the first row and last column deleted. Now, observe that $\widetilde{A} = \sum_{c} \widetilde{A}_{c}$ with $\widetilde{A}_c$ being zero except in the following blocks:
\begin{figure}[H]
\centering
\begin{tikzpicture}[baseline=(img.base)]

  \node[inner sep=0] (img) at (0,0) {
    \includegraphics[scale=0.8]{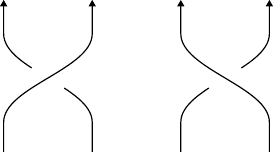}
  };

  \node at (-2, -0.8) {$i$};
  \node at (-0.8, -0.8) {$j$};
  \node at (-2.1, 0.9) {$j^{-}$};
  \node at (-0.9, 0.9) {$i^{-}$};
  \node at (0.4, -0.8) {$j$};
  \node at (1.6, -0.8) {$i$};
  \node at (0.3, 0.9) {$i^{-}$};
  \node at (1.5, 0.9) {$j^{-}$};
  \node at (-1.25, -1.5) {$s = +1$};
  \node at (1.25, -1.5) {$s = -1$};

  \node[right=0.45cm of img] (arrow) {\(\longrightarrow\)};

  \node[right=0.45cm of arrow, baseline] {
  \(
  \begin{array}{c|cc}
  \widetilde{A}_c & \text{column } i & \text{column } j \\ \hline
  \text{row } i^- & -T^{2s}  & 0  \\
  \text{row } j^- & T^{2s}-1 & -1
  \end{array}
  \)
  };
\end{tikzpicture}
\end{figure}
Map $T^2 \mapsto T$. Upon turning the knot upside down, the matrix $\widetilde{A}_c$ corresponds to the matrix $A_c$ defined in Eq. \eqref{eq:alexandermatrix}. Since the Alexander polynomial is independent of orientation (and mirror image) up to a factor $T^n$ for some $n \in \mathbb{N}$, the result follows.
\end{proof}
\begin{remark}
Although we defined the Alexander polynomial in this chapter via the matrix $A$, it is equally possible to prove Theorem~\ref{thm:alex} by defining the Alexander polynomial using the Burau representation. As a strategy, one can perform a computation analogous to that of Example~\ref{ex:firstalexander}. This approach has been carried out in~\cite{vo2017metamonoidsfoxmilnorcondition}, in the equivalent language of Gassner calculus (see Example~\ref{ex:stitching}).
  \end{remark}
  \section{Discussion}
\label{sec:discussion}
\noindent
As indicated in the introduction, the main advantage of formulating the Alexander polynomial in the language of Gaussians and contractions is that it provides a setting for introducing perturbations. To perturb the universal $R$-matrix of the XC-algebra $\mathcal{A}$ or $\mathcal{B}$, one must at least ensure that the perturbed $R$-matrix still satisfies the Yang–Baxter equation. A first approach is to add an arbitrary perturbation term to the $R$-matrix, impose the Yang–Baxter equation on this ansatz, and then solve the resulting system of equations. This “brute-force’’ method is conceptually simple, but quickly becomes cumbersome for computing higher-order perturbations of the R-matrix.

A more structured strategy exploits the representation theory of quantum groups. This is particularly natural because quantum groups carry a canonical structure of an XC-algebra \cite[Proposition 4.4]{becerra2025refinedfunctorialuniversaltangle}. One could, for example, express the Verma module in terms of Heisenberg operators, as in \cite{awata1994heisenberg}, substitute these operators into the universal $R$-matrix, and then apply the Baker--Campbell--Hausdorff formula to obtain a perturbed $R$-matrix. The zeroth-order term of the resulting expansion precisely recovers the XC-algebra $\mathcal{B}$. This observation is the main reason for working with the XC-algebra $\mathcal{B}$ rather than $\mathcal{A}$. A version of this strategy was employed in \cite{overbay2013perturbative}, though expressed in a different formalism.

There is a good reason for our insistence on using the presentation matrix of the Alexander module as given in~\eqref{eq:alexandermatrix}. Indeed, for this particular matrix, one can, in principle, obtain a closed-form expression for the perturbed Alexander invariant as formulated in \cite{barnatan2024perturbedalexanderinvariant}. From a physics perspective, one would say that, by an application of Wick’s theorem (or Theorem \ref{thm:contraction}), a perturbed Alexander invariant results from an appropriate quadratic expression in the entries of $A^{-1}$. Given an appropriate perturbation of the universal $R$-matrix in the XC-algebra $\mathcal{A}$, or $\mathcal{B}$, this follows from Theorem \ref{thm:contraction} in conjunction with Proposition \ref{pro:handlingperturbations}.

Finally, we emphasize that, despite our formulation in terms of XC-algebras, one should not lose sight of the many alternative definitions of the Alexander polynomial. In particular, descriptions that make more explicit use of topology lead to theorems that say more about the structure of a given knot. This includes, for example, the Fox-Milnor condition. To obtain similar results for perturbed Alexander invariants, it will therefore be necessary also to develop formulations that are not phrased in terms of XC-algebras, but instead in a broader language. 

\section*{Statements and Declarations}

\subsection*{Funding}
This work was supported by the Dutch Research Council (NWO) under grant number 613.009.152.

\subsection*{Conflicts of interest}
The author declares no competing interests relevant to the content of this article.

\subsection*{Availability of data and material}
Data sharing does not apply to this article as no datasets were generated or analysed during the current study.

\bibliographystyle{halpha-abbrv} 
\bibliography{bibliografia}
\appendix
\section{Useful identity to compute perturbations}
\label{sec:usefulidentity}
Throughout this paper, we have used the contraction theorem, assuming a purely Gaussian input. When we consider perturbed Gaussians, however, we must also incorporate a polynomial $P(r,s)$. When considering perturbation, it is often useful to introduce an auxiliary variable $\epsilon$. Hence, we assume $P(r,s)$ to be a series in $\epsilon$.

 By applying Proposition~\ref{thm:gencont} together with the contraction theorem, additional contractions of the form
\[
  \left\langle P \big(r+g\widetilde W , \widetilde W(s+f)\big) \right\rangle_{r,s}
\]
  must be computed. In the case of the Heisenberg algebra, the matrix $\widetilde{W}$ equals $(\mathds{1} - W)^{-1}$, with $W$ as in \eqref{eq:doubleu}. Since perturbations of the $R$-matrix depend only on the generators $x_i$ and $p_i$, the polynomial $P(r,s)$ does not depend on $\pi_i$ and $\xi_i$. As a consequence, it is only necessary to compute the bottom-left block of $(\mathds{1}-W)^{-1}$, to allow for contractions in bulk. In the following discussion, we will express this matrix in terms of the matrix $\hat{A}$ defined in Proposition~\ref{pro:contrbulk}.

Let $U_n \in \mathbb{R}^{n \times n}$ denote the upper-triangular matrix of
all ones, i.e., $(U_n)_{ij} = 1$ if $i \le j$ and $0$ otherwise. 
Let $S_n \in \mathbb{R}^{n \times n}$ denote the superdiagonal shift matrix,
defined by $(S_n)_{ij} = 1$ if $j = i+1$ and $0$ otherwise, i.e.,
\[
  S_n = \begin{pmatrix}
    0 & 1 & 0 & \cdots & 0 \\
    0 & 0 & 1 & \cdots & 0 \\
    \vdots & & \ddots & \ddots & \vdots \\
    0 & 0 & \cdots & 0 & 1 \\
    0 & 0 & \cdots & 0 & 0
  \end{pmatrix}.
\]

We note that $U_n^{-1} = \mathbb{I}_n - S_n$, which is a matrix with $1$'s on the diagonal and $-1$'s on the superdiagonal.

\begin{proposition}
\label{pro:handlingperturbations}
Let
$A \in \mathcal{A}_{n}(\bk)$. Denote by $\hat{A} \in \mathcal{A}_{n-1}(\bk)$
the matrix obtained from $A$ by deleting its first row and its last column. 
Let $M_{n}$ be as in Theorem~\ref{thm:gencont}, and define
\[
  W =
  \begin{pmatrix}
    0 & A - \mathds{1} \\
    M_{n} & 0
  \end{pmatrix}.
\]
Then the $n\times n$ bottom-left block matrix of $(\mathds{1}-W)^{-1}$ is equal to
\[
D S_n
\]
with
\begin{equation}\label{eq:shift-inverse}
D = \begin{pmatrix}
   (\mathbb{I}_{n-1}-\hat A)^{-1} & (\mathbb{I}_{n-1}-\hat A)^{-1}\hat c \\[2pt]
   0 & 1
\end{pmatrix}.
\end{equation}
where $\hat c \in \bk \llbracket h \rrbracket ^{n-1}$ denotes the last column of $A$ with its first entry
removed.
\end{proposition}

\begin{proof}

The bottom-left $n\times n$ block-matrix of $(\uno-W)^{-1}$ is of the form
\begin{equation}
\label{eq:matrix-identity}
\begin{aligned}
(\mathbb{I}_n - M_n(A-\mathbb{I}_n))^{-1}M_n
    &= (\mathbb{I}_n - (U_n - \mathbb{I}_n)(A-\mathbb{I}_n))^{-1}M_n \\
    &= (U_n(\mathbb{I}_n-A) + A)^{-1}M_n \\
    &= ((\mathbb{I}_n-A) + U_n^{-1}A)^{-1}U_n^{-1}M_n \\
    &= ((\mathbb{I}_n-A) + (\mathbb{I}_n-S_n)A)^{-1}
       (\mathbb{I}_n-S_n)M_n \\
    &= (\mathbb{I}_n-S_nA)^{-1}(\mathbb{I}_n-S_n)M_n \\
    &= (\mathbb{I}_n-S_nA)^{-1}S_n .
\end{aligned}
\end{equation}
where in the last step we used $(\mathbb{I}_n-S_n)M_n = S_n$.

To identify the first factor with the matrix $D$, let $a_1^\top,\dots,a_n^\top$ denote the rows of $A$. The matrix $S_n A$ acts on $A$ as follows:
\[
S_n A = \begin{pmatrix} a_2^\top \\ \vdots \\ a_n^\top \\ 0 \end{pmatrix},
\qquad\text{which implies}\qquad
\mathbb{I}_n - S_n A = \begin{pmatrix}
   e_1^\top - a_2^\top \\ \vdots \\ e_{n-1}^\top - a_n^\top \\ e_n^\top
\end{pmatrix}.
\]
Upon partitioning off the last row and last column, the
rows $e_i^\top - a_{i+1}^\top$ for $i=1,\dots,n-1$ give rise to the block
$\mathbb{I}_{n-1}-\hat A$ in the upper-left corner and the column $-\hat c$ in the upper-right corner. Hence, we find
\[
\mathbb{I}_n - S_n A = \begin{pmatrix} \mathbb{I}_{n-1}-\hat A & -\hat c \\[2pt] 0 & 1 \end{pmatrix}.
\]
Using the formula
\[
\begin{pmatrix} B & c \\ 0 & 1 \end{pmatrix}^{-1}
= \begin{pmatrix} B^{-1} & -B^{-1}c \\ 0 & 1 \end{pmatrix}
\]
we obtain
\begin{equation*}
(\mathbb{I}_n - S_n A)^{-1} = \begin{pmatrix}
   (\mathbb{I}_{n-1}-\hat A)^{-1} & (\mathbb{I}_{n-1}-\hat A)^{-1}\hat c \\[2pt]
   0 & 1
\end{pmatrix}. 
\end{equation*}
Substituting this expression into \eqref{eq:matrix-identity} proves the proposition.
\end{proof}

\noindent
\typeout{get arXiv to do 4 passes: Label(s) may have changed. Rerun}
\end{document}